\documentclass[11pt]{amsart}
\usepackage{amssymb}
\usepackage{amsfonts}
\usepackage{amsmath}
\usepackage[matrix,arrow,graph]{xy}

\def\thesubsection{\arabic{subsection}}
\theoremstyle{plain}
\newtheorem{thm}{Theorem}
\newtheorem{lemma}[thm]{Lemma}
\newtheorem{prop}[thm]{Proposition}

\theoremstyle{definition}

\newtheorem*{rmks}{Remarks}

\newcommand{\wt}{\widetilde}
\newcommand{\ol}{\overline}

\def\id{{\rm id}}

\newcommand{\udot}{{\scriptscriptstyle \bullet}}
\newcommand{\C}{{\mathbb C}}
\newcommand{\D}{{\mathbb D}}

\newcommand{\A}{{\mathbb A}}
\newcommand{\Z}{{\mathbb Z}}
\newcommand{\R}{{\mathbb R}}

\newcommand{\PP}{{\mathbb P}}
\newcommand{\Q}{{\mathbb Q}}
\newcommand{\F}{{\mathbb F}}
\newcommand{\HH}{{\mathbb H}}
\newcommand{\cal}{\mathcal}

\newcommand{\cL}{{\cal L}}

\newcommand{\cD}{{\cal D}}

\newcommand{\IC}{{\mathbf {IC^{\udot}}}}
\newcommand{\mb}{\mathbf}

\begin{document}
\title{Hyperbolic localization of intersection cohomology}

%\date{\today}

\author{Tom Braden}
\address{University of Massachusetts, Amherst}
\email{braden@math.umass.edu}

\begin{abstract}
For a normal variety $X$ defined over an algebraically closed field
with an action of the multiplicative group $T = G_m$, we
consider the ``hyperbolic localization'' functor $D^b(X) \to D^b(X^T)$,
which localizes using closed supports in the directions
flowing into the fixed points, and compact supports in the directions flowing
out.  We show that the hyperbolic localization of the intersection
cohomology sheaf
is a direct sum of intersection cohomology sheaves.
\end{abstract}

\maketitle
\newcommand{\sbullet}{{\scriptstyle \bullet}}
\renewcommand\thesubsection{\arabic{subsection}}

\subsection{Introduction}
Let $X$ be a normal variety over an algebraically closed field 
$k$ on which the the multiplicative group 
$T = G_m(k) = k^\times$ acts algebraically.
%, where either $k = \C$
%or the algebraic closure of a finite field $\F_q$.  
Let $F = X^T$ be the variety of fixed points, with connected components
$F_1,\dots,F_r$.   Define varieties
\[X^+_k = \{x \in X \mid \lim_{t\to 0} t\cdot x \in F_k\}, \text{ and}\]
\[X^-_k =  \{x \in X \mid \lim_{t\to \infty} t\cdot x \in F_k\}.\]
Let $X^+$ (respectively $X^-$) 
be the disjoint (disconnected) union of all the $X^+_k$ (resp. $X^-_k$),
and let $f^\pm\colon F \to X^\pm$ and
$g^\pm\colon X^\pm \to X$ restrict to the inclusion maps on each component
$F_k$, $X^\pm_k$.

If $k = \C$, let $\cD(X)$ the category of constructible complexes of
sheaves of on $X$ over a commutative ring $R$ of finite global dimension.
Otherwise, let $\cD(X)$ be either the derived category of \'etale 
sheaves on $X$ with torsion coefficients prime to the characteristic of $k$,
or one of the categories $D^b(X;\Z_l)$, $D^b(X;\Q_l)$, $D^b(X;\ol{\Q_l})$
obtained from them by standard limiting arguments.  

Define {\em hyperbolic localization} functors 
$(\mathop{\sbullet})^{!*}, 
(\mathop{\sbullet})^{*!} \colon \cD(X) \to \cD(F)$
by
\[S^{!*} := (f^+)^!(g^+)^*S,\]\[S^{*!} := (f^-)^*(g^-)^!S.\]

In general these functors are different; our first theorem gives
sufficient conditions for them to be canonically isomorphic.  
 Let $\mu\colon T \times X
\to X$ be the map defining the action.  Call an object 
$S \in \cD(X)$ ``weakly equivariant'' if $\mu^*S \cong L \boxtimes S$
for some locally constant sheaf $L$ on $T$.  Equivariant 
perverse sheaves, and more generally objects of the equivariant derived 
category of \cite{BeL} give rise to weakly equivariant objects.
Call $S$ ``$T$-constructible'' if it is 
constructible with respect to a stratification by $T$-invariant subvarieties.
Weakly equivariant objects are $T$-constructible, but not conversely.

\begin{thm} \label{iso} There is a natural morphism  $\iota_S\colon S^{*!} 
\to S^{!*}$, which is an isomorphism if $S$ is weakly equivariant, or
if $k =\C$, $R$ is an algebraically closed field, 
and $S$ is $T$-constructible.
\end{thm}

Roughly speaking, the stalks of $S^{!*}$ and $S^{*!}$ 
can be thought of as the local cohomology
of $S$ at a point of $F$  
with mixed supports --- closed supports in the directions
flowing into $F$,  and compact supports in
directions flowing away.  These two kinds of supports come from limiting
operations that do not in general commute; the content of Theorem
\ref{iso} is that they do commute for weakly equivariant objects. 

To motivate our second result, consider the case where $X$ is a
smooth projective variety over $\C$.  Then
the $F_k$ are all smooth, and the varieties 
$X^+_k$, $X^-_k$ are affine bundles over them, by \cite{B-B}.
Localizing the constant sheaf $\Q_X$ gives
\[(\Q_X)^{!*} = \bigoplus_{k=1}^r \Q_{F_k}[-n_k],\]
where $n_k = \dim_\R X^+_k - \dim_\R F_k$, and  
there is an isomorphism between $H^\udot(X;\Q)$ 
and the hypercohomology 
\[\HH^\udot((\Q_X)^{!*}) = \bigoplus_{k=1}^r H^{\sbullet - n_k}(F_k;\Q).\]

Kirwan \cite{Ki} generalized this result
to the intersection cohomology $IH^\udot(X;\Q)$ of a 
singular complex projective variety $X$:  there is a 
non-canonical isomorphism $IH^\udot(X;\Q) \cong  
\HH^\udot(F;\IC(X;\Q)^{!*})$.  Here $\IC(X)$ 
is the intersection cohomology complex of \cite{GM}.

To make Kirwan's result useful, one wants to be able to compute $\IC(X)^{!*}$. 
Our second result shows that it is as simple as one could hope.  
We assume either that $k = \C$ or else that $k = \ol{\F_q}$ and 
$X$ is obtained by extension of scalars from a variety defined over
$\F_q$.  If $k=\C$, we take $\IC(X) = \IC(X;\Q)$; if 
$k=\ol{\F_q}$, we take $\IC(X)$ to be the intersection cohomology complex
in $D^b(X;\ol{\Q_l})$.

\begin{thm} \label{main theorem} There is an isomorphism

\[\IC(X)^{!*} \cong \bigoplus_{j=1}^n \IC(Y_j,\cL_j)[d_j]\]
where $Y_1,\dots,Y_n$ are irreducible closed subvarieties of $F$, 
each $\cL_j$ is a local system on a Zariski open subset of 
the smooth locus of $Y_j$, and $d_j \in \Z$.
\end{thm}

\begin{rmks} 
(1) To see why either weak equivariance or constructibility is 
needed in Theorem \ref{iso}, take $X = \C^2$ with the action 
$t\cdot(x,y) = (tx, t^{-1}y)$, 
and let $S = i_*\Q_\C$, where $i\colon \C\to X$ is the diagonal embedding.
Then $S^{*!} = j^!S$ and $S^{!*} = j^*S$, where $j$ is the inclusion of 
the origin into $X$.
%; these are not isomorphic. 

(2) Replacing the torus action by the opposite action interchanges
$X^+$ and $X^-$, so the hyperbolic localization functors 
are sent to their Verdier duals.  Thus up to duality there is only
one hyperbolic localization functor for weakly equivariant sheaves.

(3) In \cite{GM} Goresky and MacPherson calculated local contributions to 
the Lefschetz fixed-point formula for a self-map 
(or more generally a correspondence) acting on a complex of sheaves on 
a singular space.  They assume
that the map is ``weakly hyperbolic'', which ensures that attracting and
repelling directions can be defined.  
They give a sufficient condition (Proposition 9.2) for
$\HH^\udot(S^{!*})$ and $\HH^\udot(S^{*!})$ to be isomorphic.
This condition is satisfied when the self-map 
comes from a torus action;  this can be used to show that
Theorem \ref{iso} holds when $k = \C$ and $S$ is $T$-constructible,
without the assumption on the coefficient ring $R$.

(4) %Our proof of the first statement of Theorem \ref{iso} 
%uses only the standard adjunctions and proper and smooth base change,
%so it is valid for $l$-adic sheaves, and in any setting where these 
%properties hold, such as derived categories of finite \'etale sheaves.
%The extension to $T$-constructible objects uses the complex topology in
%an essential way, however.
Theorem \ref{main theorem} is proved by using Theorem \ref{iso} to 
show that hyperbolic localization preserves purity of 
weakly equivariant mixed Hodge modules or mixed
$l$-adic sheaves.  Similar arguments were used in \cite{Gi} \cite{BrM}
in the special case when the action is completely attracting 
near a fixed point component; in this case Theorem \ref{iso} is obvious.

(5) Except for this special case, it
 is not possible to replace hyperbolic localization by ordinary
restriction in Theorem \ref{main theorem}.  
For instance, take $X$ to be a quadric
cone in $M = \C^4$ with equation $xy-zw = 0$.  Let the $T$-action be
given by $t\cdot(x,y,z,w) = (tx,t^{-1}y,z,w)$.  Then $F\cap X = X^T$ is
the union of two lines, and $\IC(X)|_F \simeq \Q_{F\cap X} \oplus \Q_0[-2]$,
which is not a sum of intersection cohomology sheaves.
A further computation with this example also shows 
that none of the functors
$(f^+)^!$, $(f^-)^*$, $(g^+)^*$, $(g^-)^!$ preserve purity individually.

\end{rmks}

Finally, we mention a useful property of hyperbolic localization, which 
is used together with Theorem \ref{main theorem} 
in \cite{BiBr} to give inequalities on Kazhdan-Lusztig polynomials.
We restrict to the case of sheaves on $X$ with coefficients in a field.
\begin{prop} Hyperbolic localization preserves local Euler characteristics:
given any point $x\in F$ and any $S\in \cD(X)$, we have 
\[\chi(S_x) = \chi((S^{!*})_x).\]
\end{prop}

Since $(f^+)^! = \D_F(f^+)^*\D_{X^+}$, this
follows from the fact that Verdier duality preserves local Euler 
characteristics.  One way to see this is to note that $\cD(X)$ is
generated as a triangulated category by (shifted) intersection cohomology 
sheaves, and then use the fact that $\D_X\IC(Y) = \IC(Y)$.  
For varieties over $\C$,
this can also be deduced from \cite{Su} or \cite[Exercise IX.12]{KS}.

\subsubsection*{Acknowledgments} The author would like to thank 
Mikhail Finkelberg, David Massey, Ivan Mirkovi\'c, Kari Vilonen,
and the referees for useful comments and corrections, 
and Nick Katz for pointing out
that Lemma \ref{special} does not hold in the \'etale setting.

\subsection{Reduction to the affine case} The following is an easy consequence
of proper base change.
\label{reduction}
\begin{lemma} \label{functorial}
The hyperbolic localization functors commute with 
pullbacks by inclusions of open $T$-invariant subvarieties 
and with pushforwards by inclusions of closed $T$-equivariant subvarieties.
\end{lemma}

%Since closed immersions are proper, the second statement comes
%from proper base-change.  If $U \subset X$ is an open subvariety,
%then $U^T = F \cap U$ and $U^+$ is an open subvariety of $X^+ \cap U$,
%giving the first statement.

To construct the natural morphism $\iota_S\colon S^{*!} \to S^{!*}$
of Theorem \ref{iso},
it is enough to construct its restriction to each component 
of $F$.  By removing all the components of $F$ but one and using 
the lemma, we can assume that 
$X$ has only one fixed point component.  In this case
$X^+ \cap X^- = F$, so the maps
$f^\pm$, $g^\pm$ form a Cartesian square.  
Let $S^+ = (g^+)_*(g^+)^*S$, and let 
$\beta\colon S \to S^+$ be the adjunction morphism.
Base change gives a natural isomorphism
\begin{eqnarray*}\label{cd}
(S^+)^{*!} & = & (f^-)^*(g^-)^!(g^+)_*(g^+)^*S \\
               & \simeq & (f^-)^*(f^-)_*(f^+)^!(g^+)^*S \simeq S^{!*}
\end{eqnarray*}
Define $\iota^{}_S$ to be the composition of (\ref{cd})  
with $\beta^{*!}\colon S^{*!} \to (S^+)^{*!}$.  
It is compatible with the pullback by the inclusion
$j\colon U \to X$ of a Zariski open subset: we have 
$j^*\iota^{}_S \simeq \iota^{}_Uj^*$.

\begin{lemma} \label{linear charts}
Any normal $T$-variety defined over an algebraically closed
field has a covering  by $T$-invariant affine open subvarieties, 
each of which is equivariantly isomorphic to a $T$-invariant closed 
subvariety of an affine space $\A^N$ on which $T$ acts linearly.
\end{lemma}

This is a consequence of results of Sumihiro \cite{Sum}.
(We remark that this is the only place where we need 
that $k$ is algebraically closed.  
Theorem 1 will still hold, for instance, when $X$ is a real variety 
for which the conclusion of Lemma \ref{linear charts} can be verified 
directly -- for instance, if $X$ is a 
$T_{\mathit max}(\R)$-invariant subvariety of 
a real flag variety $G(\R)/B(\R)$, where $G$ is a semisimple 
algebraic group, $B$ is a Borel subgroup, and 
$T_{\mathit max}$ is a maximal split torus, all
defined over $\R$).

As a result, it is enough
to prove Theorem \ref{iso} when $X = \A^N$ with a linear $T$-action.

\subsection{Comparing pushforwards and pullbacks} 
The proofs of Theorem \ref{iso} and Theorem \ref{main theorem} use the
following technical lemma.
Suppose that $T = G_m$ acts trivially on a variety $Y$, and 
that $W \to Y$ is a $T$-equivariant line bundle, provided with an
equivariant splitting $W = W_+ \oplus W_-$ 
so that all the weights of $W_+$ 
are larger than all the weights of $W_-$.  

Let $E = \PP(W) \setminus \PP(W_+)$ and  $B = \PP(W_-)$, where 
$\PP(W)$ denotes the projective bundle of lines in $W$.
Then 
$T$ acts on $E$ so that $\lim_{t\to 0} t\cdot x \in B$ for all $x\in E$.
The resulting map $p\colon E \to B$ is given by $p([w_+, w_-]) = [0,w_-]$.
Let $\phi\colon B\to Y$ and 
$i\colon B\to E$ be the projection and inclusion maps, respectively.

\begin{lemma} \label{contracting} 
There are natural isomorphisms
in $\cD(Y)$
\[ \phi_*p_*S \simeq \phi_*i^*S,\;\; 
\phi_*p_!S \simeq \phi_*i^!S\]
for any weakly equivariant object $S \in \cD(E)$.
\end{lemma}

The proof is in \S5.

If $W_-$ is one-dimensional, we get the following special
case: $\phi\colon B \to Y$ is an isomorphism, and $E$ is an
equivariant vector bundle over $B$ with only positive weights.
This case of the lemma is well-known; see \cite{Sp}, for instance.  
The general case when the
action on $B$ is non-trivial seems to be new, and is essential for the
proof of Theorem \ref{iso}.

This gives another description of hyperbolic localization.
With the notation from the introduction,
define projection maps $\pi^{\pm}\colon X^\pm \to F$ by 
$\pi^+(x) = \lim_{t\to 0} t\cdot x$ and $\pi^-(x) = \lim_{t\to \infty} 
t\cdot x$.

Adjunction gives natural morphisms
\[ (\pi^-)_*\to (f^-)^*, \;\text{and}\, (f^+)^! \to (\pi^+)_!.\]
They are isomorphisms for weakly equivariant objects,
by Lemma \ref{linear charts} and special case of 
Lemma \ref{contracting}.
Thus for any weakly equivariant 
$S \in \cD(X)$ there are natural isomorphisms 
\begin{equation}\label{aa}
S^{!*}\simeq (\pi^+)_!(g^+)^*S,\;\; S^{*!} 
\simeq (\pi^-)_*(g^-)^!S.
\end{equation}

\subsection{Proof of Theorem \ref{iso}}
First note that if Theorem \ref{iso} holds for two objects $S, S' \in \cD(X)$, 
then it holds for the cone of any morphism $S\to S'$.  The following lemma
then allows us to reduce the proof of Theorem \ref{iso}
to the case that $S$ is weakly equivariant.

\begin{lemma} \label{special} If $k = \C$ and $R$ is an algebraically 
closed field, then the 
smallest full subcategory of $\cD(X) = D^b(X;R)$ containing all 
weakly equivariant objects which is closed under taking cones is
the category of $T$-constructible objects.
\end{lemma}
\begin{proof} 
One can show that any local system on a 
$T$-invariant variety has a filtration with weakly equivariant subquotients, 
by using the Jordan decomposition of the monodromy generated by the 
action of a loop in $T$.  The lemma follows from this by induction on strata.
\qed
\end{proof}

Thus we can assume 
that $S$ is weakly equivariant.  
By the discussion of Section \ref{reduction} we can also
assume that $X = V$ is a $k$-vector space 
with a linear $T$-action.
If $V = V_+ \oplus V_- \oplus V_0$ is 
the decomposition of $V$ into subspaces where $T$ acts with
positive, negative, and zero weights, respectively, 
then we have $F = V_0$, $X^+ = V_+ \oplus V_0$, and $X^- = V_- \oplus V_0$.

Let $j$ be the inclusion of the complement of $X^+$ into $X$, and
consider the triangle containing the morphism $\beta$ from
\S2:
\[j_!j^*S \to S \stackrel{\beta}{\to} (g^+)_*(g^+)^*S \stackrel{[1]}{\to}.\]
To show that $\iota_S$ is an isomorphism it will be enough to show that 
$(j_!j^*S)^{*!}$ vanishes, since then $\beta^{*!}$ is an isomorphism.  
Thus we can assume that $S|_{X^+} = 0$.

Let $\ol{X} = V_0 \times (\PP(V_+\oplus V_- \oplus k) \setminus \PP(
V_-))$; define a $T$-action by letting $T$ act trivially on the summand $k$.  
There is an equivariant open embedding $b\colon X\to \ol{X}$ 
given by $b(v_+, v_-, v_0) = v_0 \times [v_+, v_-, 1]$.

The closure of $b(X^+)$ in $\ol{X}$ is $V_0 \times \PP(V_+\oplus  k)$;
denote it by $\ol{X^+}$, and let $\ol{g^+}$ be the inclusion into $\ol{X}$.
Let $\pi^-$, $\ol{\pi^+}$, $\ol{\pi}$ be the projections of $X^-$, 
$\ol{X^+}$, and $\ol{X}$, respectively onto $V_0$.

Note that $b(X^-)$ is closed in $\ol{X}$; 
let $h$ be the inclusion of the complement.  There is
a distinguished triangle:
\begin{equation}\label{triangle}
b_*(g^-)_*(g^-)^!S \to b_!S \to h_*h^*b_!S \stackrel{[1]}{\to}\; 
\end{equation}
Applying $\ol{\pi}_*$, the left term becomes  
\[\ol{\pi}_*b_*(g^-)_*(g^-)^!S = (\pi^-)_*(g^-)^!S \simeq S^{*!},\]
using (\ref{aa}).  Theorem \ref{iso} will follow if we can show that 
$\ol{\pi}_*$ annihilates the second and 
third terms of (\ref{triangle}).

For the second
term we see that $E = \ol{X}$, $B = \ol{X^+}$ are of the form
given in Lemma \ref{contracting}, so we 
get $\ol{\pi}_*b_!S = (\ol{\pi^+})_*(\ol{f^+})^*b_!S$.  This
vanishes, since $S|_{X^+} = 0$ implies $(b_!S)|_{\ol{X^+}}=0$.  

For the third term we use Lemma \ref{contracting} again, this time with
$E = \ol{X} \setminus X$,  $B = \ol{X^+} \setminus X^+$.  This is a
vector bundle with only negative weights, so using the special case
of Lemma \ref{contracting} (with the opposite $T$-action), we get
$\ol{\pi}_*h_*h^*b_!S \simeq \ol{\pi}_*a^*b_!S = 0$, where
$a$ is the inclusion of $\ol{X^+} \setminus X^+$ into $\ol{X}$.

\subsection{Proof of Theorem 2} \label{mhms}
\newcommand{\bL}{{\mb L}}

We recall the properties of mixed Hodge modules and mixed
$l$-adic sheaves that we will need to prove Theorem 2; 
more information can be found in \cite{BBD}, \cite{S}.

If $k = \C$, we let $\cD_m(X)$ be Saito's category of mixed Hodge
modules.  If $k = \ol{\F_q}$ and $X$ arises by extension of scalars
from a variety $X_0$ defined over $\F_q$, we let $\cD_m(X)$ be
the category $D^b_{\mathit mixed}(X_0;\ol{\Q_l})$ of mixed $\ol{\Q_l}$ sheaves
on $X_0$.  These two situations share the following formal properties.

There is a triangulated functor $r\colon\cD_m(X) \to \cD(X)$
with the property that $r(S) = 0$ if and only if $S = 0$; it follows that
any morphism $\alpha$ in $\cD_m(X)$ is an 
isomorphism if and only if $r(\alpha)$ is.

For any $w\in \Z$, there are full subcategories $\cD_{\le w}(X)$
and $\cD_{\ge w}(X)$ of $\cD_m(X)$.  Objects of 
$\cD_{\le w}(X) \cap \cD_{\ge w}(X)$ are called {\em pure
of weight $w$}.  If $Y$ is a subvariety of $X$ and $w\in \Z$, 
the intersection cohomology sheaf 
$\IC(Y) \in \cD(X)$ with constant coefficients
extended by $0$ to $X$ is isomorphic to $r(L)$ 
for an object $L \in \cD_m(X)$ which is pure of weight $w$.
Conversely, if $S$ is pure, then $r(S)$ is a direct sum of intersection
cohomology sheaves with possibly non-constant coefficients.

Given an algebraic map $f\colon X\to Y$, the 
functors $f^*$, $f^!$, $f_*$, $f_!$ lift to functors between 
$\cD_m(X)$ and $\cD_m(Y)$ which satisfy the same adjunction properties
as their non-mixed versions.  We have
\[f_*\cD_{\ge w}(X) \subset \cD_{\ge w}(Y),\; f^!\cD_{\ge w}(Y) \subset \cD_{\ge w}(X)\]
\[f_!\cD_{\le w}(X) \subset \cD_{\le w}(Y),\; f^*\cD_{\le w}(Y) \subset \cD_{\le w}(X)\]
We say an object $S \in \cD_m(X)$ is weakly equivariant or $T$-constructible
if the same property holds for $r(S)$.

Using these properties, we can now deduce Theorem \ref{main theorem} from 
Theorem \ref{iso}.  
The isomorphisms (\ref{aa}) and the transformation $\iota^{}_S$
are defined by means of the 
standard adjunctions, so
they lift to isomorphisms on
weakly equivariant mixed sheaves.
Thus, for any $w \in \Z$ we have 
\[\cD_{\le w}(X)^{!*} = 
(\pi^+)_!(g^+)^*\cD_{\le w}(X) \subset \cD_{\le w}(F),\text{ and}\]
\[\cD_{\ge w}(X)^{!*} = \cD_{\ge w}(X)^{*!} = 
(\pi^-)_*(g^-)^!\cD_{\ge w}(X) \subset \cD_{\ge w}(F).\] 
We have proved:

\begin{thm}
The hyperbolic localization functors preserve purity of weakly equivariant
mixed sheaves.
\end{thm}
Theorem \ref{main theorem} follows.

\subsection{Proof of Lemma \ref{contracting}}
 We will prove the first isomorphism of  Lemma \ref{contracting}; 
the second follows by dualizing.  The following proof is a generalization
of an argument of Springer \cite{Sp}.
   
Applying $\phi_* p_*$ to the adjunction map $S \to i_*i^*S$ 
gives a natural map $\phi_*p_*S \to \phi_*i^*S$.
If it is an isomorphism for two objects in a distinguished triangle,
it will also be an isomorphism for the third.  If $j$ is the
inclusion of $E \setminus i(B)$ into $E$, we have a distinguished
triangle:
\[j_!j^*S \to S \to i_*i^*S \stackrel{[1]}{\to}\; \]
Since Lemma \ref{contracting} 
is clearly true for $i_*i^*S$, it is enough to 
show it when $S=j_!j^*S$.  In this case $i^*S = 0$, 
and we must prove that $\phi_*p_*S = 0$.

Let  $\Gamma_\mu$ be the graph of the action map
$\mu\colon T \times E \to E$, and let 
$\Gamma$ be its closure in $\A^1 \times E \times E$.
Let $q_1$, $q_2\colon \Gamma\to \A^1 \times E$ 
be the restrictions of the projections $\pi_{12}$ and $\pi_{13}$.  
An easy computation in coordinates shows that if  
$\wt{\Gamma}$ is the closure of $\Gamma_\mu$ in 
$\A^1 \times E \times \PP(W)$, then $\wt{\Gamma} \setminus \Gamma_\mu
\subset \{0\}\times E \times i(B)$.  As a result, $\wt{\Gamma} = \Gamma$, 
and the map $q_1$ is proper.  %This is the only place where we use the
%specific construction of the spaces $E$ and $B$. 

Let $K$ denote the coefficient ring of our sheaves.
We will write $K_Z$ to denote the constant rank $1$ local system
on a variety $Z$, placed in degree $0$.
Define a map 
$\psi = \id_{\A^1} \times (\phi p)\colon \A^1 \times E \to \A^1 \times Y$.
Adjunction gives a morphism
\begin{align}\label{key morphism}
\psi_*(K_{\A^1} \boxtimes S) & \to  
\psi_*q_{2*}q_2^*(K_{\A^1} \boxtimes S) \\ \nonumber & \simeq 
\psi_*q_{1!}q_2^*(K_{\A^1} \boxtimes S),
\end{align}
using the properness of $q_1$ and the fact that  $\psi q_1 = \psi q_2$.
The object on the left is isomorphic to $K_{\A^1} \boxtimes \phi_*p_*S$,
by smooth base change.  Since $q_2$ restricts to an isomorphism 
over $T \times E$, we see that $\alpha$ is an isomorphism
on $T \times Y$.
 
Let $\wt S$ be the right hand side of (\ref{key morphism}).
We will show that $\wt S|_{\{0\} \times Y} = 0$.  Since
$\wt{S}|_{T\times Y} \simeq K_T \boxtimes \phi_*p_*S$, this
implies that $\wt S \simeq b_!K_{T} \boxtimes \phi_*p_*S$, where
$b\colon T \to \A^1$ is the inclusion. 

Composing the morphism (\ref{key morphism}) with the
morphism $b_!K_{T} \boxtimes \phi_*p_*S \to K_{\A^1} \boxtimes \phi_*p_*S$
obtained from the adjunction morphism $b_!K_T \to K_{\A^1}$
gives an endomorphism of $K_{\A^1} \boxtimes \phi_*p_*S$
which is an isomorphism over $T \times Y$ but which restricts to $0$ on
$\{0\} \times Y$.  Since this is only possible if $\wt S = 0$, we 
must have $\phi_*p_*S = 0$, as desired.

Finally, to show that 
$\wt S|_{\{0\} \times Y} = 0$,  let $E_0 = E\setminus i(B)$, and
consider the maps
\[\xymatrix{
T \times E_0 \ar[r]^{\hat \jmath}\ar[d]^\pi & \Gamma \ar[d]^{\pi'q_2}\\
E_0 \ar[r]^j & E}
\]
where $j$ is the inclusion,  
$\hat \jmath(t,x) = (t, t^{-1}\cdot x, x)$ and $\pi, \pi'$ are
the projections from $T \times E_0$ and $\A^1\times E$
onto the second factors. The square commutes,
and is in fact Cartesian, since 
$(\pi' q_2)^{-1}(E_0) \subset \Gamma_\mu$.

Let $S_0 = S|_{E_0}$.  Then $S \simeq j_!S_0$, so
\[\wt S \simeq 
\psi_*q_{1!}q_2^*(\pi')^*j_!S_0 \simeq 
\psi_*q_{1!}\hat{\jmath}_!\pi^*S_0.\]
Since $S$, and thus $S_0$, is weakly equivariant, we have
\[(q_{1!}\hat{\jmath}_!\pi^*S_0)|_{T \times E_0} \simeq L \boxtimes S_0\]
for some local system $L$ on $T$.
This implies $q_{1!}\hat{\jmath}_!\pi^*S_0 \simeq b_!L \boxtimes 
j_!S_0 = b_!L \boxtimes S$, since the stalks of this sheaf vanish outside 
of $T \times E_0$.  Applying $\psi_*$ now gives
\[\wt{S} \simeq b_!L \boxtimes \phi_*p_*S,\]
which implies $\wt S|_{\{0\} \times Y} = 0$, as required.

\end{document}